\documentclass{amsart}

\usepackage{latexsym}
\usepackage{amssymb}
\usepackage{amsmath}
\usepackage{amsthm}
\usepackage{fancyhdr}
\usepackage{setspace}
\usepackage[colorlinks=true,linkcolor=blue,citecolor=blue,urlcolor=blue]{hyperref}

\newcommand{\RR}{{\mathbb R}}

\newcommand{\tr}{\mathrm{tr}}

\newcommand{\lgra}{\longrightarrow}

\newtheorem{example}{Exemples}[section]
\newtheorem{thm}{Theorem}[section]
\newtheorem*{thrm}{Theorem}
\newtheorem{lemma}[thm]{Lemma}

\newtheorem{cor}[thm]{Corollary}
\newtheorem{definition}[thm]{Definition}
\newtheorem{notation}[thm]{Notation}
\newtheorem{exabout:ample}[thm]{Example}

\theoremstyle{definition}
\newtheorem{remark}[thm]{Remark}
\newtheorem{remarks}[thm]{Remarks}

\newcommand{\beqt}{\begin{equation}}  \newcommand{\eeqt}{\end{equation}}
\newcommand{\bal}{\begin{align}}      \newcommand{\eal}{\end{align}}
\newcommand{\ba}{\begin{array}}      \newcommand{\ea}{\end{array}}
\newcommand{\bc}{\begin{center}}     \newcommand{\ec}{\end{center}}
\newcommand{\be}{\begin{enumerate}}  \newcommand{\ee}{\end{enumerate}}
\newcommand{\beq}{\begin{eqnarray}}  \newcommand{\eeq}{\end{eqnarray}}
\newcommand{\beQ}{\begin{eqnarray*}} \newcommand{\eeQ}{\end{eqnarray*}}
\newcommand{\bi}{\begin{itemize}}    \newcommand{\ei}{\end{itemize}}
\newcommand{\bt}{\begin{tabular}}    \newcommand{\et}{\end{tabular}}

\renewcommand{\beqt}{\begin{equation}\begin{alignedat}{2}}
\renewcommand{\eeqt}{\end{alignedat}\end{equation}}

\begin{document}
\author{Julien Roth and Julian Scheuer}

\title{Pinching of the first eigenvalue for second order operators on hypersurfaces of the Euclidean space}
\date{\today}

\address{Julien Roth, Laboratoire d'Analyse et de Math\'ematiques Appliqu\'ees, UPEM-UPEC, CNRS, F-77454 Marne-la-Vall\'ee, France}
\email{julien.roth@u-pem.fr}

\address{Julian Scheuer, Albert-Ludwigs-Universit{\"a}t, Mathematisches Institut, Eckerstr.~1,
79104 Freiburg, Germany}
\email{julian.scheuer@math.uni-freiburg.de}

\begin{abstract}
We prove stability results associated with upper bounds for the first eigenvalue of certain second order differential operators of divergence-type on hypersurfaces of the Euclidean space. We deduce some applications to $r$-stability as well as to almost-Einstein hypersurfaces.\end{abstract}

\maketitle

\setcounter{tocdepth}{1}

\tableofcontents
\section{Introduction}
There is a wide literature concerning estimates of the eigenvalues of the Laplacian (and more general divergence-type second order operators) on submanifolds of spaceforms. The first upper bound for the first eigenvalue of the Laplacian on compact hypersurfaces of $\RR^{n+1}$ was obtained by Bleecker and Weiner \cite{BleeckerWeiner:/1976} who showed 
\beqt\label{inegBW}
\lambda_1\leqslant \frac{1}{V(M)}\int_M|B|^2dv_g,
\eeqt
where $B$ is the second fundamental form of the hypersurface $M$ and $V(M)$ its volume. After that, Reilly \cite{Reilly:/1977} improved this upper bound by getting the norm of the mean curvature instead of the second fundamental form. Precisely, he proved the following estimate:
\beqt\label{inegReH}
\lambda_1\leqslant \frac{n}{V(M)}\int_M|H|^2dv_g.
\eeqt
For these two inequalities the limiting case is attained if and only if the hypersurface is a hypersphere. Later on, Heintze \cite{Heintze:04/1988} extended Reilly's upper bound to hypersurfaces of compact ambient spaces and, after a partial result due to Heintze, El Soufi and Ilias \cite{El-SoufiIlias:/1992} proved an analogue in the hyperbolic space, also cf.~\cite{AlencarCarmoMarques:/2001}. In fact, in \cite{Reilly:/1977} Reilly proved a sequence of  upper bounds involving higher order mean curvatures, for which \eqref{inegReH} is just the particular case  $r=0$:
\beqt\label{inegReHr}
\lambda_1\left(\int_MH_rdv_g\right)^2\leqslant nV(M)\int_MH_{r+1}^2dv_g.
\eeqt
As the previous ones the inequalities are sharp and the limiting cases are also characterised by geodesic hyperspheres. Note that the inequalities in \eqref{inegReHr} are not in general better than \eqref{inegReH}; they are not comparable to the other.

The Laplacian is not the only fundamental operator for which extrinsic estimates have been proved. In particular comparable upper bounds have been established  for some divergence-type second order elliptic operators as the operators $L_r$ associated with the higher order mean curvatures $H_r$ and even for more general elliptic divergence-free operators, cf.~\cite{Grosjean:/2004}.  Notations and basic facts about higher order mean curvatures will be given in Section \ref{sec2}.  

The aim of the present article is first to study the stability of the limiting case for a sequence of optimal upper bounds for a larger class of second order operators including the Laplacian and the operators $L_r$. In \cite{Roth:/2015a} the first author proved the  the following theorem.
\begin{thrm}[\cite{Roth:/2015a}]
Let $(M^n,g)$ be a connected, oriented and closed Riemannian manifold isometrically immersed into $\RR^{n+1}$. Assume that $M$ is endowed with two symmetric and divergence-free $(1,1)$-tensors $S$ and $T$. Assume in addition that $T$ is positive definite. Then the first positive eigenvalue of the operator $L_T$ satisfies
\beqt\lambda_1(L_T)\left(\int_M\tr(S)dv_{g}\right)^2\leqslant \left(\int_M\tr(T)dv_g\right)\left(\int_M|H_S|^2dv_g\right).
\eeqt
Moreover, if $H_S$ does not vanish identically and equality occurs, then $M$ is a hypersphere and $\tr(S)$ is constant.
\end{thrm}
The operator $L_T$ and also $H_{S}$ will be defined in Section \ref{sec2}. Since the equality is achieved only for hyperspheres, it is natural to study the stability of the limiting case. We give the following quantitative result about the proximity to hyperspheres for hypersurfaces almost satisfying the equality case. Namely, under the condition
\beqt\label{PinchLambda}
\lambda_1(L_T)\left(\int_M\tr(S)dv_{g}\right)^2\geqslant (1-\varepsilon)\left(\int_M\tr(T)dv_g\right)\|H_S\|^2_{2p}V(M),
\eeqt
for $p>1$, we prove the following result.
\begin{thm}\label{thm1}Let $n\geqslant 2$ and $(M^n,g)$ be a connected, oriented and closed Riemannian manifold isometrically immersed into $\RR^{n+1}$ by $X$. Assume that $M$ is endowed with two symmetric and divergence-free $(1,1)$-tensors $S$ and $T$. Assume in addition that $T$ is positive definite, $S^{ij}B_{ij}$ is not identically zero and that for some $q>n$ there holds $V(M)\|B\|_{q}^{n}\leqslant A.$ Let $p>1$. Then there exists $\varepsilon_{0}=\varepsilon_{0}(n,p,q,A)>0$, such that if \eqref{PinchLambda} holds with $\varepsilon<\varepsilon_{0}$, then   $M$ is $\varepsilon^{\alpha}$-close, diffeomorphic and $\varepsilon^{\alpha}$-almost isometric to a sphere $S(\bar{X},r)$, i.e.
there exists $r>0$, a constant $C$ depending on $n$, $p$, $q$ and $A$, and $\alpha=\alpha(n,q)$, such that

\beqt\label{thm1-1}
\||X-\bar{X}|-r\|_{\infty}\leq Cr\varepsilon^{\alpha},
\eeqt
where $\bar{X}=\int_{M}X$ is the center of mass of $X(M)$,
and for a natural diffeomorphism 
\beqt
F\colon (M,d_{1})\rightarrow \left(S(r),d_{2}\right)
\eeqt
we have
\beqt
|d_{2}(F(x_{1}),F(x_{2}))-d_{1}(x_{1},x_{2})|\leqslant Cr\varepsilon^{\alpha}\quad \forall x_{1},x_{2}\in M.
\eeqt
  Moreover, $M$ is embedded and $X(M)$ is a starshaped hypersurface.
\end{thm}

\begin{remark}
In order to provide \eqref{thm1-1}, it is not necessary that $V(M)\|B\|_{q}^{n}\leqslant A$, but only $V(M)\|H\|^{n}_{q}\leqslant A.$
\end{remark}

In the sections \ref{SecStability} and \ref{SecEinstein} we derive some applications of Theorem \ref{thm1}. The first of those concerns the $r$-stability. Namely we prove that if a hypersurface with constant $r$-th mean curvature is almost stable in a suitable sense, then it is a geodesic sphere (see Theorem \ref{thm2}). Finally, the last section of the paper will be devoted to almost-Einstein hypersurfaces, where we improve previous closeness results of the first author \cite{Roth:05/2008} and Hu, Xu and Zhao \cite{HuXuZhao:/2015}. 

Earlier results on the eigenvalue pinching problem can be found in \cite{ColboisGrosjean:/2007}, \cite{Roth:05/2008} for the Euclidean case, \cite{GrosjeanRoth:/2012} for other ambient spaces and \cite{HuXuZhao:/2015} for an improvement of parameters compared to \cite{Roth:05/2008}.

\section{Preliminaries}\label{sec2}
\subsection{Generalized Hsiung-Minkowski formula}
Let $(M^n,g)$ be a connected, oriented and closed Riemannian manifold isometrically immersed into $\RR^{n+1}$. We denote by $X$ its vector position, by $\nu$ its unit normal vector and we consider a divergence-free symmetric $(1,1)$-tensor $T$ on $M$. We associate to $T$ the second order differential operator $L_T$ defined by
$$L_Tu:=-\mathrm{div}(T\nabla u),$$
for any $C^2$-function $u$ on $M$. We also associate to $T$ the following generalised mean curvature vectors:
\beqt\label{defHT}
H_T=T^{ij}X_{ij},
\eeqt
where $X_{ij}$ is the second covariant derivative of the immersion vector $X$. Then we have
\beqt\label{LTX}
L_TX=-H_T
\eeqt
and easily deduce the following identity:
\beqt\label{LTX2}
\frac{1}{2}L_T|X|^2=-\langle X,H_T\rangle-\tr(T).
\eeqt
After integrating this gives the so-called generalized Hsiung-Minkowski formula
\beqt\label{genHM}
\int_M\Big(\langle X,H_T\rangle+\tr(T)\Big)dv_g=0.
\eeqt
The classical Hsiung-Minkowski formula is in fact this formula for the particular case where $T$ is the identity and more generally, where $T$ is the tensor $T_r$ associated to the $r$-th mean curvature. In the next paragraph we will recall the definition of $T_r$ and the associated operator $L_r$, which play a crucial role for the $r$-stability as we will see in Section \ref{SecStability}.

\subsection{Higher order mean curvatures}
The higher order mean curvatures are extrinsic quantities defined from the second fundamental form and generalising the notion of mean curvature. Up to a normalisation constant the mean curvature $H$ is the trace of the second fundamental form $B$:
\beqt
H=\frac{1}{n}\tr(B).
\eeqt
In other words the mean curvature is  
\beqt
H=\frac{1}{n}S_1(\kappa_1,\dots,\kappa_n),
\eeqt
where $S_1$ is the first elementary symmetric polynomial and $\kappa_1,\dots,\kappa_n$ are the principal curvatures. Higher order mean curvatures are defined in a similar way for $r\in\{1,\dots,n\}$ by
\beqt 
H_r=\frac{1}{\binom{n}{r} }S_r(\kappa_1,\cdots,\kappa_n),
\eeqt
where $S_r$ is the $r$-th elementary symmetric polynomial, that is for any $n$-tuple $(x_1,\cdots,x_n)$, 
\beqt
S_r(x_1,\dots,x_n)=\sum_{1\leqslant i_1<\dots<i_r \leqslant n}x_{i_1}\cdots x_{i_r}.
\eeqt
 By convention we set $H_0=1$ and $H_{n+1}=0$. Finally, for convenience we also set $H_{-1}=-\langle X,\nu\rangle$.

To each $H_r$ we associate a symmetric $(2,0)$-tensor, which is in coordinates given by
\beqt
T_r=(T^{ij}_{r})=\left(\frac{\partial S_{r+1}}{\partial B_{ij}}\right),
\eeqt
where $S_{r+1}$ is now understood to depend on the second fundamental form and the metric. The relation between these two notions can be found in \cite[Sec.~2.1]{Gerhardt:/2006} for example. These tensors $T_r$ are divergence-free (see \cite{Grosjean:09/2002} for instance) and satisfy the following relations:
\beqt 
\tr(T_r)=c(r)H_r\quad \text{and}\quad H_{T_r}=-c(r)H_{r+1}\nu,
\eeqt
where $c(r)=(n-r)\binom{n}{r}$ and $H_{T_r}$ is given by the relation (\ref{defHT}).\\
We finish this section by giving some classical inequalities between the $H_r$ which are well-known. First, for any $r\in\{0,\cdots,n-2\}$, 
\beqt\label{inegHr1}
H_rH_{r+2}\leqslant H_{r+1}^2,
\eeqt 
with equality at umbilical points, cf.~\cite[p.~104]{HardyLittlewoodPolya:/1952}. Moreover, cf.~ \cite{BarbosaColares:06/1997}, if $H_{r+1}>0,$ then $H_s>0$ for any $s\in\{0,\cdots,r\}$ and 
\beqt\label{inegHr2}
H_{r+1}^{\frac{1}{r+1}}\leqslant H_r^{\frac{1}{r}}\leqslant \cdots\leqslant H_2^{\frac12}\leqslant H.
\eeqt
Combining \eqref{inegHr1} and \eqref{inegHr2}, we get that  if $H_{r+1}>0$, then
\beqt\label{inegHr3}
H_{r+2}\leqslant HH_{r+1},
\eeqt with equality at umbilical points.

\subsection{Extrinsic Sobolev inequality}
We finish this section of preliminaries by recalling the extrinsic Sobolev inequality proved by Michael and Simon \cite{MichaelSimon:05/1973} for submanifolds of the Euclidean space and by Hoffman and Spruck \cite{HoffmannSpruck:11/1974}  for any ambient space. This inequality will play a crucial role  in the iteration process to obtain $L^{\infty}$-estimates in the proof of Theorem \ref{thm1}. For any positive $C^1$-function on a hypersurface $M$ of $\RR^{n+1}$ we have
\begin{equation}\label{simon}
\left(\int_Mf^{\frac{n}{n-1}}dv_g\right)^{\frac{n-1}{n}}\leqslant K_n\int_M\left(|\nabla f|+|H|f\right)dv_g,
\end{equation}
where $K_n$ is a positive constant depending only on the dimension $n$ and $H$ is the mean curvature of $M$. An immediate consequence of this inequality is that we can bound the volume of $M$ in terms of the mean curvature from below. Indeed we have $V(M)^{\frac{n-1}{n}}\leqslant K_n\|H\|_{\infty} V(M)$, obtained by taking  $f\equiv1$ in \eqref{simon}, which gives the following lower bound for the volume:
\beqt\label{minvolume}
V(M)\geqslant\frac{1}{(K_n\|H\|_{\infty})^n}.
\eeqt
Finally we fix our convention for the $L^p$-norms. For $p\in [1,\infty)$ and a function $f$ defined on $M$ we set
\beqt
\|f\|_p=\left(\frac{1}{V(M)}\int_M|f|^pdv_g\right)^{\frac1p}.
\eeqt

\section{A general upper bound for \texorpdfstring{$\lambda_1(L_T)$}{}}\label{sec3}
We first recall the general upper bound for the first eigenvalue of the operator $L_T$ in terms of a symmetric  divergence-free $(1,1)$-tensor $S$.

\begin{thm}\label{thm0}{\cite{Roth:/2015a}}
Let $(M^n,g)$ be a connected, oriented and closed Riemannian manifold isometrically immersed into $\RR^{n+1}$. Assume that $M$ is endowed with two symmetric and divergence-free $(1,1)$-tensors $S$ and $T$. Assume in addition that $T$ is positive definite. Then the first positive eigenvalue of the operator $L_T$ satisfies
\beqt\label{ineggene}
\lambda_1(L_T)\left(\int_M\tr(S)dv_{g}\right)^2\leqslant \left(\int_M\tr(T)dv_g\right)\left(\int_M|H_S|^2dv_g\right).
\eeqt
Moreover, if $H_S$ does not vanish identically and equality occurs, then $M$ is a hypersphere and $\tr(S)$ is constant.
\end{thm}
{\bf Proof:} First we recall that $L_TX=-H_T$. Without loss of generality let 
\beqt\int_{M}X=0.
\eeqt
 Thus we may use the coordinates $X^{\alpha}$ as test functions in the Rayleigh quotient to obtain
\beqt
\lambda_1(L_T)\int_M|X|^2dv_g&\leqslant-\int_M\langle X,H_T\rangle dv_g\\
&=\int_M\tr(T)dv_g,\label{quotient}
\eeqt
where we have used the generalised Hsiung-Minkowski formula \eqref{genHM} for the tensor $T$. Hence, we deduce that
\beqt
\left(\int_M\tr(T)dv_g\right)\left(\int_M|H_S|^2dv_g\right)&\geqslant\lambda_1(L_T)\left(\int_M|X|^2dv_g\right)\left(\int_M|H_S|^2dv_g\right)\\
&\geqslant\lambda_1(L_T)\left(\int_M|X\|H_S|dv_g\right)^2\\
&\geqslant\lambda_1(L_T)\left(\int_M\langle X,H_S\rangle dv_g\right)^2\\
&=\lambda_1(L_T)\left(\int_M\tr(S) dv_g\right)^2,
\eeqt
where we have used the Cauchy-Schwarz inequality and the generalised Hsiung-Minkowski formula for the tensor $S$.\\
Let's study the limiting case. If $H_S$ does not vanish identically, the integral of $\tr(S)$ is not zero in the equality case. Hence we have equality in the Cauchy-Schwarz inequality, that is, $X$ and $H_S$ are colinear. In particular the position vector is normal. We deduce easily that $M$ is a sphere. Indeed we have
\beqt
\frac{\partial}{\partial \xi^{i}}\langle X,X\rangle=2\langle X_{i},X\rangle=0.
\eeqt
 Thus the norm of $X$ is constant and so $M$ is contained in a sphere. Since $M$ has no boundary, $M$ is the entire sphere. Moreover $M$ is totally umbilic and so by definition of $H_S$ we have $H_S=-k\tr(S)\nu$. Finally $\tr(S)=k^{-1}|H_S|$ is constant, since by equality in the Cauchy-Schwarz inequality $H_S$ and  $X$ are proportional. 
\hfill $\square$

\begin{remarks}
\begin{enumerate}
\item Note that for the equality case there is no equivalence.  Indeed, if equality holds, then $M$ is necessarily a geodesic sphere, but there is no reason that equality occurs if $M$ is a geodesic sphere.
\item If $T$ is the identity, that is $L_T$ is the Laplacian, we have an alternative proof of the limiting case. If equality holds, then the coordinates are eigenvectors of the Laplacian and $M$ is a geodesic sphere by Takahashi's Theorem \cite{Takahashi:/1966}.
\end{enumerate}
\end{remarks}

The classical inequalities are particular cases of this general one. First of all, for $T= {\rm Id}$ we get the Reilly inequalities \eqref{inegReH} and \eqref{inegReHr} cited above by taking $S={\rm Id}$ or $S=T_r$. For $T=T_r$, $0\leqslant r\leqslant n-1$, we recover the following inequalities for the  first eigenvalue of the operator $L_r$ obtained by Alias and Malacarne \cite{AliasMalacarne:/2004}:
 \beqt\label{inegAM} \lambda_1(L_r)\left(\int_MH_sdv_g\right)^2\leqslant c(r)\int_MH_rdv_g\int_MH_{s+1}^2dv_g,
 \eeqt
 with $S=T_s$. The special case $r=s$ gives a well known inequality proved by  Alencar, do Carmo and Rosenberg \cite{AlencarCarmoRosenberg:11/1993}:
 \beqt\label{inegACR} \lambda_1(L_r)\int_MH_rdv_g\leqslant c(r)\int_MH_{r+1}^2dv_g.\eeqt
Another particular case of our general inequality is the case $S=T$. We get the following Corollary.

 \begin{cor}
 Let $(M^n,g)$ be a closed, connected and oriented Riemannian manifold isometrically immersed into $\RR^{n+1}$. Assume that $M$ is endowed with a divergence-free and positive definite symmetric $(1,1)$-tensor $T$.  Then, the first positive eigenvalue of the operator $L_T$ satisfies
 \beqt
 \lambda_1(L_T)\left(\int_M\tr(T)dv_{g}\right)\leqslant \int_M|H_T|^2dv_g.
 \eeqt
Moreover, if equality holds, then $M$ is a geodesic sphere and $\tr(T)$ is constant.
 \end{cor}
 For $T=\mathrm{Id}$ we get the Reilly inequality \eqref{inegReH} and for $T=T_r$ we have equality \eqref{inegACR}.

\section{A pinching result for \texorpdfstring{$\lambda_1(L_T)$}{}}\label{sec4}
In this section we consider the pinching problem associated with \eqref{ineggene}, which is the stability of its equality case. In other words, if equality is almost achieved, $M$ is close to a sphere in a sense to be made precise. For technical reasons we will consider a less sharp but nevertheless optimal inequality which admits the same equality case:
\beqt\label{ineqthm0Lp}
\lambda_1(L_T)\left(\int_M\tr(S)dv_{g}\right)^2\leqslant \left(\int_M\tr(T)dv_g\right)\|H_S\|^2_{2p}V(M),
\eeqt
for $p>1$, obtained from \eqref{ineggene} by H\"older inequality. We introduce the following pinching condition for $0<\varepsilon<1$:
\begin{equation}
\lambda_1(L_T)\left(\int_M\tr(S)dv_{g}\right)^2\geqslant (1-\varepsilon)\left(\int_M\tr(T)dv_g\right)\|H_S\|^2_{2p}V(M)\label{pincLT} \tag{$\Lambda_{p,\varepsilon}$}.
\end{equation}

\subsection{\texorpdfstring{$L^2$}{L^2}-estimates}
In this section we prove some lemmata giving proximity in an $L^2$-sense between the hypersurface $M$ and a geodesic sphere of appropriate radius under the assumptions of Theorem \ref{thm1}. 

\begin{lemma}\label{XL2}
If \eqref{pincLT} is satisfied and $\bar{X}=0$, then
\beqt\label{XL21}
\frac{\int_M\tr(T)dv_g}{\lambda_1(L_T)V(M)}(1-\varepsilon)^2\leqslant \|X\|^2_2\leqslant \frac{\int_M\tr(T)dv_g}{\lambda_1(L_T)V(M)}
\eeqt
and 
\beqt\label{XL22} 
\frac{\left(\frac{1}{V(M)}\int_{M}\tr(S)dv_{g}\right)^{2}}{\|H_{S}\|_{2p}^{2}}(1-\varepsilon)^{2}\leqslant \|X\|^2_2\leqslant\frac{1}{1-\varepsilon}\frac{\left(\frac{1}{V(M)}\int_{M}\tr(S)dv_{g}\right)^{2}}{\|H_{S}\|_{2p}^{2}}.
\eeqt
\end{lemma}

{\bf Proof:}
We start from \eqref{quotient} obtained by injecting the coordinate functions into the Rayleigh quotient. This gives the upper bound for $\|X\|_2^2$ in \eqref{XL21}. For the lower bound we start from \eqref{quotient} again. We have
\beqt
&\lambda_1(L_T)\int_M|X|^2dv_g\left(\int_M\tr(S)dv_g\right)^4\\
\leqslant&\int_M\tr(T)dv_g\left(\int_M\tr(S)dv_g\right)^4\\
=&\int_M\tr(T)dv_g\left(\int_M\langle H_S,X\rangle dv_g\right)^4\\
\leqslant&\int_M\tr(T)dv_g\left(\int_M|H_S|^2dv_g\right)^2\left(\int_M |X|^2 dv_g\right)^2.
\eeqt
Hence we get
\beqt\label{XL2a}
\|X\|_2^2\geqslant\frac{\lambda_1(L_T)\left(\int_M\tr(S)dv_g\right)^4}{\left(\int_M\tr(T)dv_g\right)\|H_S\|_{2p}^4V(M)^3}.
\eeqt
Now we use \eqref{pincLT}, which gives
\beqt
\lambda_1(L_T)^2\geqslant \frac{\left(\int_M\tr(T)dv_g\right)^2\|H_S\|_{2p}^4V(M)^2}{\left(\int_M\tr(S)dv_g\right)^4}(1-\varepsilon)^2
\eeqt
and, together with \eqref{XL2a}, yields the desired result,
\beqt
\frac{\int_M\tr(T)dv_g}{\lambda_1(L_T)V(M)}(1-\varepsilon)^2\leqslant \|X\|^2_2.\eeqt
The upper estimate of \eqref{XL22} follows from inserting \eqref{pincLT} into \eqref{XL21}, the lower one from inserting \eqref{ineqthm0Lp} into \eqref{XL21}.
\hfill$\square$

Now we state a second lemma which gives an $L^2$-estimate for the tangential part of the position vector under the pinching condition.

\begin{lemma}\label{TangL2}
If \eqref{pincLT} holds and $\bar{X}=0$, then
\beqt
\|X^T\|^2_2\leqslant\varepsilon\|X\|^2_2.
\eeqt
\end{lemma}

{\bf Proof:} As in the previous lemma we start from \eqref{quotient}. We have
\beqt
&\lambda_1(L_T)\int_M|X|^2dv_g\left(\int_M\tr(S)dv_g\right)^2\\
\leqslant&\int_M\tr(T)dv_g\left(\int_M\tr(S)dv_g\right)^2\\
=&\int_M\tr(T)dv_g\left(\int_M\langle H_S,X\rangle dv_g\right)^2\\
\leqslant&\int_M\tr(T)dv_g\left(\int_M|H_S|^2dv_g\right)\left(\int_M \langle X,\nu\rangle^2 dv_g\right)\\
\leqslant&\int_M\tr(T)dv_g\left(\int_M \langle X,\nu\rangle^2 dv_g\right)\|H_S\|_{2p}^2V(M).
\eeqt
From this we deduce that
\beqt
\|X^T\|_2^2&=\frac{1}{V(M)}\int_M(|X|^2-\langle X,\nu\rangle^2)dv_g\\
&\leqslant\|X\|_2^2-\dfrac{\lambda_1(L_T)\|X\|_2^2\left(\int_M\tr(S)dv_g\right)^2}{\left(\int_M\tr(T)dv_g\right)\|H_S\|_{2p}^2V(M)}\\
&=\|X\|_2^2\left(1-\dfrac{\lambda_1(L_T)\left(\int_M\tr(S)dv_g\right)^2}{\left(\int_M\tr(T)dv_g\right)\|H_S\|_{2p}^2V(M)}\right)\\
&\leqslant\varepsilon\|X\|_2^2,
\eeqt
by using \eqref{pincLT} for the last line.
\hfill$\square$

\subsection{From $L^2$ to $L^{\infty}$}
Now we will give a sequence of lemmas, based on an iteration process, which allow us to control the $L^{\infty}$-norm of some functions by their $L^2$-norm. Note that this iteration process does not depend on the pinching condition. We have the following lemma for the norm of the position vector. The proof can be found in \cite[Lemma~5]{HuXuZhao:/2015}. This lemma is an improvement of a similar lemma given in \cite{GrosjeanRoth:/2012} and \cite{Roth:12/2013}. 

\begin{lemma}\cite[Lemma~5]{HuXuZhao:/2015}\label{SphClose}
Let $q>n$ be a real number. There exists a constant $\Gamma(n,q)>0$, so that for any isometrically immersed, compact submanifold $M^{n}$ of $\RR^{n+1}$ we have
\beqt
\||X-\bar{X}|-\|X-\bar{X}\|_{2}\|_{\infty}\leqslant\Gamma \left(V(M)\|H\|_q^n\right)^{\frac{\gamma}{n}}\|X-\bar{X}\|_2\left(1-\frac{\|X-\bar{X}\|_1}{\|X-\bar{X}\|_2}\right)^{\frac{1}{2(\gamma+1)}},
\eeqt
where $\gamma=\frac{nq}{2(q-n)}$.
\end{lemma}

The following $L^{\infty}$ estimate of the tangential part is a straightforward adaption of the proof of \cite[Lemma~6]{HuXuZhao:/2015}. Let us include it with the necessary modifications for completeness.

\begin{lemma}\label{TangSmall}
Let $q>n$ and $X\colon M\hookrightarrow \RR^{n+1}$ be the immersion of a closed hypersurface. Then there exists a constant $C=C(n,q),$ such that 
\beqt
\|X^{T}\|_{\infty}\leqslant C\left(V(M)\|B\|_{q}^{n}\|X\|_{\infty}\right)^{\frac{\gamma}{\gamma+1}}\|X^{T}\|_{2}^{\frac{1}{\gamma+1}}.
\eeqt
\end{lemma}

{\bf Proof:}
Set $\psi=|X^{T}|.$ Due to 
\beqt
\psi^{2}=|X|^{2}-\langle X,\nu\rangle^{2}
\eeqt
we obtain
\beqt
|d\psi|\leqslant 1+n|B|\|X\|_{\infty}.
\eeqt
Applying \eqref{simon} with $f=\psi^{2\alpha}$, $\alpha\geq 1$, we obtain
\beqt
\|\psi\|_{\frac{2\alpha n}{n-1}}^{2\alpha}\leqslant KV(M)^{\frac 1n}\left(2\alpha (1+n\|B\|_{q}\|X\|_{\infty})+\|H\|_{q}\|\psi\|_{\infty}\right)\|\psi\|^{2\alpha-1}_{\frac{(2\alpha-1)q}{q-1}}.
\eeqt 
From here we follow the proof of \cite[Lemma~5]{HuXuZhao:/2015}. Setting
\beqt
\mu=\frac{n(q-1)}{(n-1)q},\quad a_{0}=\frac{2q}{q-1},\quad a_{p+1}=\mu a_{p}+\frac{n}{n-1},
\eeqt
then we obtain, using
\beqt
\alpha=a_{p}\frac{q-1}{2q}+\frac 12,
\eeqt
that
\beqt
&\left(\frac{\|\psi\|_{a_{p+1}}}{\|\psi\|_{\infty}}\right)^{\frac{a_{p+1}}{\mu^{p+1}}}\\
\leqslant &\left(C(n)V(M)^{\frac 1n}\left(a_{p+1}\frac{1+n\|B\|_{q}\|X\|_{\infty}}{\|\psi\|_{\infty}}+\|H\|_{q}\right)\right)^{\frac{n}{\mu^{p+1}(n-1)}}\left(\frac{\|\psi\|_{a_{p}}}{\|\psi\|_{\infty}}\right)^{\frac{a_{p}}{\mu^{p}}}\\
\leqslant&\left(C(n)V(M)^{\frac 1n}a_{p+1}\frac{\|B\|_{q}\|X\|_{\infty}}{\|\psi\|_{\infty}}\right)^{\frac{n}{\mu^{p+1}(n-1)}}\left(\frac{\|\psi\|_{a_{p}}}{\|\psi\|_{\infty}}\right)^{\frac{a_{p}}{\mu^{p}}},
\eeqt
where we also used
\beqt
1\leq \|X\|_{2}\|H\|_{2}\leq C(n)\|X\|_{\infty}\|B\|_{q},\quad \|\psi\|_{\infty}\leq \|X\|_{\infty}.
\eeqt
As in the proof of \cite[Lemma~5]{HuXuZhao:/2015} we obtain
\beqt
\|\psi\|_{\infty}\leqslant C(q,n)\left(V(M)^{\frac 1n}\frac{\|B\|_{q}\|X\|_{\infty}}{\|\psi\|_{\infty}}\right)^{\gamma}\|\psi\|_{2},
\eeqt
where
\beqt
\gamma=\frac{nq}{2(q-n)}.
\eeqt
Rearranging terms gives the desired estimate.
\hfill$\square$

\subsection{Proof of Theorem \ref{thm1}}

In order to prove Theorem \ref{thm1}, we examine the consequences from the $L^{\infty}$-estimates under the pinching condition \eqref{pincLT}. The proofs closely follow the argumentation in \cite[Sec.~3]{HuXuZhao:/2015}.

\begin{lemma}\label{L1L2}
Under the condition \eqref{pincLT} the immersion $X$ satisfies
\beqt
1-\frac{\|X-\bar{X}\|_{1}}{\|X-\bar{X}\|_{2}}\leq C(p)\varepsilon.
\eeqt
\end{lemma}

{\bf Proof:}
Assume by translation that $\bar{X}=0$ and use Lemma \ref{XL2} and the Hsiung-Minkowski formula  to deduce
\beqt
\|H_{S}\|_{2p}\|X\|_{2}&\leqslant\frac{\int_{M}|\langle X,H_{S}\rangle|}{V(M)\sqrt{1-\varepsilon}}\\
				&\leqslant\frac{1}{\sqrt{1-\varepsilon}}\|H_{S}\|_{2p}\|X\|_{\frac{2p}{2p-1}}\\
				&\leqslant \frac{1}{\sqrt{1-\varepsilon}}\|H_{S}\|_{2p}\|X\|_{1}^{1-\frac 1p}\|X\|_{2}^{\frac 1p}.
\eeqt
We obtain
\beqt
\|X\|_{2}\leqslant (1-\varepsilon)^{-\frac{p}{2(p-1)}}\|X\|_{1}
\eeqt
and hence
\beqt
1-\frac{\|X\|_{1}}{\|X\|_{2}}\leqslant 1-(1-\varepsilon)^{\frac{p}{2(p-1)}}\leqslant C(p)\varepsilon
\eeqt
with a constant $C=C(p).$
\hfill$\square$

Let us combine these results to get final pinching estimates for $|X|$ and $|X^{T}|$. First we prove \eqref{thm1-1}.

\begin{cor}\label{PinchX}
Let $q>n$ and suppose $V(M)\|H\|^{n}_{q}\leqslant A.$
Define
\beqt
r=\frac{\frac{1}{V(M)}\left|\int_{M}\mathrm{tr}(S)\right|}{\|H_{S}\|_{2p}}
\eeqt
and let $S$ and $T$ satisfy the conditions of Theorem \ref{thm1}.
Then under the condition \eqref{pincLT} with $p>1$ and $\varepsilon<\frac 12$ there exists a constant
$C=C(n,p,q,A)$, such that
\beqt\label{PinchX1}
\||X-\bar{X}|-r\|_{\infty}\leq Cr\varepsilon^{\frac{1}{2(\gamma+1)}}.
\eeqt
\end{cor}

{\bf Proof:}
Again let $\bar{X}=0.$ Then we obtain
\beqt
\||X|-r\|_{\infty}\leqslant \||X|-\|X\|_{2}\|_{\infty}+|\|X\|_{2}-r|.
\eeqt
The lemmata \ref{XL2}, \ref{SphClose} and \ref{L1L2} give the result.

\hfill$\square$

\begin{cor}\label{PinchTang}
Let $q>n$, $V(M)\|B\|^{n}_{q}\leqslant A$ and $\bar{X}=0.$ Let $S$ and $T$ satisfy the conditions of Theorem \ref{thm1}. Then under the condition \eqref{pincLT} with $p>1$ and $\varepsilon<\frac 12$ there exists a constant $C=C(n,p,q,A)$, such that
\beqt\label{PinchTang1}
\|X^{T}\|_{\infty}\leqslant Cr\varepsilon^{\frac{1}{2(\gamma+1)}}.
\eeqt
\end{cor}

{\bf{Proof:}}
Due to Lemma \ref{TangL2} and Lemma \ref{TangSmall} we have
\beqt
\|X^{T}\|_{\infty}\leqslant C\|X\|_{\infty}^{\frac{\gamma}{\gamma+1}}\|X^{T}\|_{2}^{\frac{1}{\gamma+1}}\leqslant C\|X\|_{\infty}^{\frac{\gamma}{\gamma+1}}\|X\|_{2}^{\frac{1}{\gamma+1}}\varepsilon^{\frac{1}{2(\gamma+1)}}.
\eeqt
Since
\beqt
\|X\|_{\infty}\leqslant \||X|-r\|_{\infty}+r
\eeqt
and 
\beqt \|X\|_{2}\leqslant Cr,
\eeqt
we obtain the result from Corollary \ref{PinchX}.
\hfill$\square$

\begin{remark}\label{ProofMainThm}
Now the proof of the rest of Theorem \ref{thm1} proceeds literally as the corresponding proof of \cite[Thm.~2]{HuXuZhao:/2015}. The starting points are the equations \eqref{PinchX1} and \eqref{PinchTang1}, which correspond to the equations (10) and (12) in \cite{HuXuZhao:/2015}. Our radius $r$ corresponds to the quantity $\frac{\|H_{k-1}\|_{1}}{\|H_{k}\|_{p}}$ in \cite{HuXuZhao:/2015}. Note that this is where we have to restrict to some small $\varepsilon_{0}>0$ as claimed in Theorem \ref{thm1}.
\end{remark}

\section{Application to \texorpdfstring{$r$}{r}-stability}\label{SecStability}
In this section, we are interested in the stability of constant mean curvature hypersurfaces and, more generally, in $r$-stability. For this we introduce the $r$-area functionals
\beqt
\mathcal{A}_r=\left(\int_MS_rdv_g\right), \quad r\in\{0,\dots,n-1\}.
\eeqt
Now we consider a variation of the immersion $X_{0}$. Precisely,  let $\delta>0$ and 
\beqt 
X\colon (-\delta,\delta)\times M\lgra\RR^{n+1},
\eeqt
 such that for all $t\in (-\delta,\delta)$, $X_t:=X(t,\cdot)$ is an immersion of $M$ into $\RR^{n+1}$ and $X(0,\cdot)=X_{0}$. The precise assumptions on $M$ will be specified in Theorem \ref{thm2} below. We denote by $S_r(t)$ the corresponding curvature functions, by $\mathcal{A}_r(t)$ the $r$-area of  $X_t$ and finally we set 
 \beqt f_t=\left\langle \frac{dX}{dt},\nu_t \right\rangle,\eeqt
  where $\nu_t$ is the unit normal to $M$ induced by $X_t$ and $g_t$ is the induced metric on $M$.
  
Note that 
\beqt\frac{d}{dt}\left(\int_MS_r(t)dv_{g_t}\right)=-(r+1)\int_M fS_{r+1}dv_{g_{t}},\eeqt
cf.~\cite[Thm.~B]{Reilly:/1973}.
We also consider the volume functional
\beqt
V(t)=\int_{[0,t)\times M}X^*dv.
\eeqt
It is easy to see, cf. \cite[Lemma~2.1]{BarbosaCarmoEschenburg:/1988}, that $V$ satisfies
\beqt
 V'(t)=\int_Mf_tdv_{g_t}
\eeqt
and so $X$ preserves the volume if and only if $\int_Mf_{t}dv_{g_t}=0$ for all $t$. Moreover, according to \cite[Lemma~2.2]{BarbosaCarmoEschenburg:/1988}, for any function $f_0:M\rightarrow\RR$ such that $\int_Mf_0dv_g=0$, there exists a variation of $X_{0}$ preserving the volume and with  normal part given by $f_0$. Hence, in order to study variations with constant volume, we can equivalently consider normal parts with vanishing integral. This has the advantage that they can be used as test functions in the Rayleigh quotient for the first eigenvalue of the Laplacian or, more generally, for the operator $L_r$. As we will see, these operators play a crucial role in the study of $r$-stable hypersurfaces. Indeed, a standard argument as in \cite[Prop.~2.7]{BarbosaCarmo:/1984} shows that $X_{0}$ is a critical point for the functional $\mathcal{A}_r$ for variations with constant volume if and only if $H_{r+1}$ is constant. For such a critical point, Reilly \cite{Reilly:/1973} has computed the second variation of $\mathcal{A}_r$, also compare \cite[equ.~(2.2)]{AlencarCarmoRosenberg:11/1993}:
\beqt\label{var2Ar}
\mathcal{A}_r''(0)&=(r+1)\int_Mf_{0}L_rf_{0}+\left[c(r+1)H_{r+2}-n\frac{c(r)}{r+1}HH_{r+1}\right]f_{0}^2dv_g\\
			&=:J_{r}f_{0},
\eeqt
where $f_{0}$ is the normal part of a variation with constant volume at $t=0$ and 
\beqt
c(r)=(n-r)\binom{n}{r}.
\eeqt
 We say that a hypersurface $M$ with constant $H_{r+1}$ is $r$-stable if $J_{r}f$ is non-negative for all functions of vanishing integral. A classic result of Barbosa and do Carmo \cite{BarbosaCarmo:/1984} for $r=0$, Alencar, do Carmo and Colares \cite{AlencarCarmoColares:/1993} for $r=1$ and Alencar, do Carmo and Rosenberg \cite{AlencarCarmoRosenberg:11/1993} for all $r$ says that a compact hypersurface of the Euclidean space with constant $H_{r+1}$ is $r$-stable if and only if it is a geodesic sphere.  This result has been generalised later by Barbosa and Colares \cite{BarbosaColares:06/1997} for hypersurfaces of the hyperbolic space and the half-sphere. 

Here, using the pinching result for the first eigenvalue of $L_r$ proven above, we will show that this characterisation of spheres of the Euclidean space remains true if we only assume that $J_{r}f$ is almost positive, which means
\beqt
 \label{almost r-stable}J_{r}f\geqslant -\varepsilon \int_{M}f^{2}H_{r+1}^{\frac{r+2}{r+1}}\quad\forall f\in C^{\infty}(M)\colon \int_{M}f=0,
 \eeqt where $\varepsilon$ is a sufficiently small positive constant. Note that the term $H_{r+1}^{\frac{r+2}{r+1}}$ is present in order to have a condition which is invariant under any homothety of $\RR^{n+1}$. Precisely, we prove the following result.
 
\begin{thm}\label{thm2} Let $(M^n,g)$ be a connected, oriented and closed manifold isometrically immersed into $\RR^{n+1}$ with constant $H_{r+1}>0$. Let $A>0$, $q>n$ and assume that $V(M)\|B\|_q^{n}\leqslant A$. Then there exists $\varepsilon_{0}>0$ depending on $n$, $q$ and $A$ such that if $M$ is almost $r$-stable in the sense of \eqref{almost r-stable} with $\varepsilon<\varepsilon_{0}$,
then $M$ is a sphere of radius $H_{r+1}^{-\frac{1}{r+1}}$.

\end{thm}
{\bf Proof:} 
 Taking $f$ as an eigenfunction of $L_r$ associated with the first positive eigenvalue, we get from the almost-stability that
 \beqt
 0&\leqslant \int_{M}\left((r+1)\lambda_{1}+(r+1)c(r+1)H_{r+2}-nc(r)H_{1}H_{r+1}+\varepsilon H_{r+1}^{\frac{r+2}{r+1}}\right)f^{2}\\
 	&\leqslant \int_{M}\left((r+1)\lambda_{1}+\left((r+1)c(r+1)-nc(r)\right)H_{1}H_{r+1}+\varepsilon H_{r+1}^{\frac{r+2}{r+1}}\right)f^{2}\\
	&\leqslant\int_{M}\left((r+1)\lambda_{1}-((r+1)c(r)-\varepsilon)H_{r+1}^{\frac{r+2}{r+1}}\right)f^{2}\\
	&=(r+1)\int_{M}\left(\lambda_{1}-c(r)\left(1-\frac{\varepsilon}{(r+1)c(r)}\right)H_{r+1}^{\frac{r+2}{r+1}}\right) f^{2},
	 \eeqt
 where we used \eqref{inegHr2}, \eqref{inegHr3} and also
 \beqt
 c(r+1)-n\frac{c(r)}{r+1}=-c(r).
 \eeqt
Hence we find
\beqt
\lambda_{1}-c(r)\left(1-\frac{\varepsilon}{(r+1)c(r)}\right)H_{r+1}^{\frac{r+2}{r+1}}\geqslant 0.
\eeqt
 
 Up to a minor change this is the pinching condition of Theorem \ref{thm1}. Indeed, for $p>1$ we have
 \beqt
 \lambda_1(L_r)&\geqslant\left(1-\frac{\varepsilon}{(r+1)c(r)}\right)c(r)H_{r+1}^{\frac{r+2}{r+1}}\\
 			&=\left(1-\frac{\varepsilon}{(r+1)c(r)}\right)c(r)\frac{H^{2}_{r+1}}{H_{r+1}^{\frac{r}{r+1}}}\\
			&=\left(1-\frac{\varepsilon}{(r+1)c(r)}\right)c(r)\frac{\|H_{r+1}\|_{2p}^{2}V(M)}{\int_{M}H_{r+1}^{\frac{r}{r+1}}}\\
			&\geqslant \left(1-\frac{\varepsilon}{(r+1)c(r)}\right)c(r)\frac{\|H_{r+1}\|_{2p}^{2}V(M)}{\int_{M}H_{r}},
  \eeqt
 where we used the fact that $H_{r+1}$ is constant to make integrals appear and \eqref{inegHr2} for the last line. Now we have exactly the pinching condition for $\lambda_1(L_r)$ with $S=T=T_{r}$ and we can apply Theorem \ref{thm1} to conclude that in particular $M$ is embedded. In this argument also note that $T_{r}$ is elliptic since $H_{r+1}$ is positive. By the Alexandrov theorem for $H_{r+1}$ proved by Ros, cf.~\cite{Ros:/1987}, we get that $M$ is a sphere.
\hfill$\square$

\section{An almost-Einstein type estimate}\label{SecEinstein}

We finish this paper by applying the eigenvalue pinching result to deduce an explicit spherical closeness estimate of almost-Einstein type for hypersurfaces, namely an estimate of the form
\beqt\label{AlmEinstein1}d_{\mathcal{H}}(M,S(\bar{X},r))\leqslant \frac{cr}{\bar{R}^{\alpha}}\left\|\mathrm{Ric}-\frac{\bar{R}}{n}g\right\|_{p}^{\alpha},
\eeqt
whenever the right-hand side is small. We assume that the average of the scalar curvature $R$ over $M$ is positive, $\bar{R}>0$, and that $R\geq 0$. We will prove this estimate for $p>\max(2,\frac n2)$ and $c$ and $\alpha$ will be suitable constants depending on $n$, $p$ and an integral bound on $H$. In \cite[Thm.~3]{HuXuZhao:/2015} a similar application of the eigenvalue pinching was given, however here the authors derived the almost-isometry to the sphere with the help of the pinching
\beqt
\|\mathrm{Ric}-\lambda g\|_{\infty}<\varepsilon.
\eeqt
Here we want to improve this application, which is by the way also possible with the help of their eigenvalue pinching result: we want to relax the norm of the pinching quantity to $L^{p}$ instead of $L^{\infty}$. Besides the improvement of the estimate itself this has the advantage that we also obtain an almost-Schur type estimate right away, due to the well-known almost-Schur estimate by De Lellis and Topping, cf.~\cite{De-LellisTopping:/2012}. For manifolds of nonnegative Ricci curvature and $n\geqslant 3$ they provide the estimate
\beqt
\int_{M}\left|\mathrm{Ric}-\frac{\bar{R}}{n}g\right|^{2}\leqslant\frac{n^{2}}{(n-2)^{2}}\int_{M}\left|\mathrm{Ric}-\frac{R}{n}g\right|^{2}.
\eeqt
Similar estimates, relaxing the assumption on nonnegative Ricci curvature in different directions, were obtained by Cheng, cf.~\cite[Thm.~1.2]{Cheng:/2013} and by Ge and Wang in \cite{GeWang:03/2012}.
Instead of $L^{2}$ we need an $L^{p}$-estimate, but as was pointed out in \cite{De-LellisTopping:/2012}, their proof easily adapts provided one has a Calderon-Zygmund type inequality. For the sake of completeness, we will provide the proof of this $L^{p}$-estimate. For a class of manifolds which allow for a Calderon-Zygmund estimate see the paper \cite{Hiroshima:/1996}.

To achieve the relaxation from the $L^{\infty}$- to an $L^{p}$-pinching condition, we imitate the proof of a similar result by the second author concerning an almost-umbilical type estimate, compare the proof of \cite[Thm.~4.1]{Scheuer:08/2015}. A similar result due to the first author under slightly different assumptions was achieved in \cite{Roth:02/2013}.

Let us first recall De Lellis' and Topping's argument from \cite{De-LellisTopping:/2012} how to obtain the $L^{p}$-version of their almost-Schur lemma.

\begin{lemma}
Let $n\geq 3$ and $(M^{n},g)$ be a closed Riemannian manifold, which admits an estimate of the form
\beqt
\|\nabla^{2}u\|_{q}\leqslant C\|\Delta u\|_{q},\quad q>1,
\eeqt  
for all smooth functions $u$, where $C$ is a constant only depending on $n$, $q$ and possibly some fixed geometric quantities of $M$. Then for $p\geqslant 2$ there holds
\beqt
\left\|\mathrm{Ric}-\frac{\bar{R}}{n}g\right\|_{p}\leqslant c\left\|\mathrm{Ric}-\frac{R}{n}g\right\|_{p},
\eeqt
where $c=c(C,n,p)$.
\end{lemma}

{\bf{Proof}:} Compare \cite[Sec.~2]{De-LellisTopping:/2012}.
Let $f\in C^{2,p-2}(M)$ be the solution of
\beqt
\begin{cases}\Delta f=|R-\bar{R}|^{p-2}(R-\bar R)\\
	\int_{M}f=0. \end{cases}
\eeqt
Then, with exactly the same computation as in \cite[equ.~(2.3)]{De-LellisTopping:/2012}, we obtain
\beqt
\frac{1}{V(M)}\int_{M}|R-\bar{R}|^{p}&=\frac{1}{V(M)}\int_{M}(R-\bar{R})\Delta f\\
				&\leqslant \frac{2n}{n-2}\|\mathring{\mathrm{Ric}}\|_{p}\|\nabla^{2}f\|_{\frac{p}{p-1}},
\eeqt
and hence, due to the Calderon-Zygmund estimate,
\beqt
\frac{1}{V(M)}\int_{M}|R-\bar{R}|^{p}&\leqslant c\|\mathring{\mathrm{Ric}}\|_{p}\||R-\bar{R}|^{p-1}\|_{\frac{p}{p-1}}.
\eeqt
This yields the desired estimate, also using
\beqt
\left|\mathrm{Ric}-\frac{\bar R}{n}g\right|^{2}=|\mathring{\mathrm{Ric}}|^{2}+\frac{1}{n}(R-\bar R)^{2}.
\eeqt
\hfill$\square$

Let us come to the proof of \eqref{AlmEinstein1}. Here we will need to estimate the first Laplace eigenvalue in terms of an $L^{p}$-Ricci bound. We have the following estimate due to Aubry, originally proved in \cite[Prop.~1.5]{Aubry:/2007}, maybe more accessible in the version of \cite[Thm.~1.6]{Aubry:/2009}. It says that for $p>n/2$,
a complete Riemannian manifold $(M^{n},g)$ with
\beqt
\frac{1}{V(M)}\int_{M}(\underline{\mathrm{Ric}}-(n-1))_{-}^{p}<\frac{1}{C(p,n)}
\eeqt
is compact and satisfies
\beqt
\lambda_{1}\geqslant n\left(1-C(n,p)\left(\frac{1}{V(M)}\int_{M}(\underline{\mathrm{Ric}}-(n-1))_{-}^{p}\right)^{\frac 1p}\right),
\eeqt
where $\underline{\mathrm{Ric}}$ denotes the smallest eigenvalue of the Ricci tensor and $x_{-}=\max(0,-x).$

\begin{thm}\label{AlmEinstein}
Let $n\geq 2,$ $q>n$, $p>\max(2,\frac n2)$ and $(M^{n},g)$ be a closed, connected and oriented Riemannian manifold with $\bar{R}>0$ and $R\geqslant 0$, isometrically immersed into $\RR^{n+1}.$ Suppose $V(M)\|H\|^{n}_{q}\leqslant A$. Then there exists a constant $\varepsilon_{0}(n,p,q,A)>0$, such that whenever there holds
\beqt
\left\|\mathrm{Ric}-\frac{\bar R}{n}g\right\|_{p}\leqslant \varepsilon_{0}\bar{R},
\eeqt 
then we also have
\beqt
d_{\mathcal{H}}(M,S(\bar X,r))\leq \frac{cr}{\bar{R}^{\alpha}}\left\|\mathrm{Ric}-\frac{\bar R}{n}g\right\|_{p}^{\alpha},
\eeqt
where $\alpha=\alpha(n,q)$ and $c=c(n,p,q,A).$ If furthermore $V(M)\|B\|_{q}^{n}\leqslant A,$ then $M$ is embedded as a starshaped hypersurface and $\varepsilon^{\alpha}$-almost-isometric in the sense of Theorem \ref{thm1}.
\end{thm}

{\bf{Proof}:}
We want to apply Theorem \ref{thm1}. We apply this theorem with the tensors $T=\mathrm{Id}$ and $S=T_{1}.$ Hence we have to provide the estimate
\beqt
\lambda_{1}\|H\|_{1}^{2}\geqslant (1-\varepsilon)n\|H_{2}\|_{2s}^{2},
\eeqt
where $\lambda_{1}$ now is the first eigenvalue of $L_{T}=-\Delta$ and $s=\frac{p}{2}$.
Using a simple rescaling argument we obtain a scaled version of Aubry's result \cite[Thm.~1.6]{Aubry:/2009}, namely that
\beqt
\left(\frac{1}{V(M)}\int_{M}\left(\underline{\mathrm{Ric}}-\frac{\bar R}{n}\right)^{2s}_{-}\right)^{\frac{1}{2s}}<\frac{\bar R}{C(n,s)}
\eeqt
implies
\beqt\label{AlmEinsteinA}
\lambda_{1}&\geqslant \frac{\bar R}{n-1}\left(1-\frac{C(s,n)}{\bar R}\left(\frac{1}{V(M)}\int_{M}\left(\underline{\mathrm{Ric}}-\frac{\bar R}{n}\right)^{2s}_{-}\right)^{\frac{1}{2s}}\right)\\
		&\geqslant \frac{\bar R}{n-1}-C(s,n)\left\|\mathrm{Ric}-\frac{\bar R}{n}g\right\|_{2s}.
\eeqt
Let us estimate $\|H\|_{1}^{2}.$
\beqt\label{AlmEinsteinB}
\|H\|_{1}^{2}&\geqslant\left(\frac{1}{V(M)}\int_{M}H^{\frac 12}_{2}\right)^{2}\\
		&=\frac{1}{n(n-1)}\left(\frac{1}{V(M)}\int_{M}R^{\frac 12}\right)^{2}\\
		&= \frac{\bar R}{n(n-1)}+\frac{1}{n(n-1)}\left(\|R^{\frac 12}\|_{1}^{2}-\|\bar{R}^{\frac 12}\|_{1}^{2}\right)\\
		&\geqslant \frac{\bar R}{n(n-1)}-\frac{1}{n(n-1)}(\|R^{\frac 12}\|_{1}+\bar{R}^{\frac 12})\|R^{\frac 12}-\bar{R}^{\frac 12}\|_{1}\\
		&\geqslant \frac{\bar{R}}{n(n-1)}-\frac{1}{n(n-1)\sqrt{\bar R}}\|R-\bar{R}\|_{1}(\|R^{\frac 12}\|_{1}+\bar{R}^{\frac 12})\\
		&\geqslant \frac{\bar{R}}{n(n-1)}-c_{n}\left\|\mathrm{Ric}-\frac{\bar R}{n}g\right\|_{1}.
\eeqt
We also have to estimate $\|H_{2}\|_{2s}^{2}.$
\beqt\label{AlmEinsteinC}
\|H_{2}\|_{2s}^{2}&=\frac{1}{n^{2}(n-1)^{2}}\left(\frac{1}{V(M)}\int_{M}R^{2s}\right)^{\frac 1s}\\
			&\leqslant \frac{\bar{R}^{2}}{n^{2}(n-1)^{2}}+\frac{1}{n^{2}(n-1)^{2}}\left|\|R\|_{2s}^{2}-\|\bar{R}\|_{2s}^{2}\right|\\
			&\leqslant \frac{\bar{R}^{2}}{n^{2}(n-1)^{2}}+\frac{\|R\|_{2s}+\bar{R}}{n^{2}(n-1)^{2}}\|R-\bar{R}\|_{2s}\\
			&\leqslant \frac{\bar{R}^{2}}{n^{2}(n-1)^{2}}+c_{n}(\|R\|_{2s}+\bar{R})\left\|\mathrm{Ric}-\frac{\bar{R}}{n}g\right\|_{2s}.			
\eeqt
Bringing together \eqref{AlmEinsteinA}, \eqref{AlmEinsteinB} and \eqref{AlmEinsteinC}, we obtain, also noting that the right-hand side of \eqref{AlmEinsteinA} is still positive when $\varepsilon_{0}>0$ is small,
\beqt
\lambda_{1}\|H\|_{1}^{2}&\geqslant \frac{\bar{R}^{2}}{n(n-1)^{2}}-c_{n}\bar{R}\left\|\mathrm{Ric}-\frac{\bar{R}}{n}g\right\|_{1}-c(n,s)\bar{R}\left\|\mathrm{Ric}-\frac{\bar{R}}{n}g\right\|_{2s}\\
			&\geqslant n\|H_{2}\|_{2s}^{2}-c(\|R\|_{2s}+\bar{R})\left\|\mathrm{Ric}-\frac{\bar{R}}{n}g\right\|_{2s}\\
			&=\left(1-c\frac{\|R\|_{2s}+\bar{R}}{\|R\|^{2}_{2s}}\left\|\mathrm{Ric}-\frac{\bar{R}}{n}g\right\|_{2s}\right)n\|H_{2}\|_{2s}^{2}\\
			&\geqslant\left(1-\frac{c}{\bar{R}}\left\|\mathrm{Ric}-\frac{\bar{R}}{n}g\right\|_{2s}\right)n\|H_{2}\|_{2s}^{2}.
\eeqt

Thus when $\varepsilon_{0}$ is chosen small enough, we may apply Theorem \ref{thm1} to conclude the result.

\hfill$\square$

\begin{remark}
Let us note again, as already mentioned in the introduction to this section, that Theorem \ref{AlmEinstein} enables us to obtain an almost-Schur type estimate, whenever we obtain an estimate in the sense of De Lellis and Topping \cite{De-LellisTopping:/2012}. 
\end{remark}

\bibliographystyle{hamsplain}
\bibliography{bibliography}

\end{document}